\documentclass[12pt]{amsart}
\usepackage{amssymb}
\usepackage{amsfonts}
\usepackage{latexsym}
\usepackage{amscd}

\usepackage[all]{xy}

\input xy
\xyoption{all}

\numberwithin{equation}{section}

\newcommand{\Z}{{\mathbb{Z}}}

\newtheorem{lem}{Lemma}[section]
\newtheorem{corol}[lem]{Corollary}
\newtheorem{theor}[lem]{Theorem}
\newtheorem{prop}[lem]{Proposition}
\newtheorem{rema}[lem]{Remark}

\newtheorem{defis}[lem]{Definitions}
\newtheorem{exem}[lem]{Example}

\begin{document}
\title[The socle of a Leavitt path algebra]{The socle of a Leavitt path algebra}

\author{G. Aranda Pino}
\address{Centre de Recerca Matem\`{a}tica, Apartat 50, E-08193, Bellaterra (Barcelona), Spain.}
\email{garanda@crm.es}

\author{D. Mart\'\i n Barquero}
\address{Departamento de Matem\'atica Aplicada, Universidad de M\'alaga, 29071 M\'alaga, Spain.}
\email{dmartin@uma.es}

\author{C. Mart\'\i n Gonz\'alez} \author{M. Siles Molina}
\address{Departamento de \'Algebra, Geometr\'{\i}a y Topolog\'{\i}a, Universidad de M\'alaga, 29071 M\'alaga, Spain.}
\email{candido@apncs.cie.uma.es} \email{mercedes@agt.cie.uma.es}


\begin{abstract} In this paper we characterize the minimal left ideals of a Leavitt path algebra as those ones which
are isomorphic to principal left ideals generated by line point vertices, that is, by vertices whose trees do not
contain neither bifurcations nor closed paths. Moreover, we show that the socle of a Leavitt path algebra is the
two-sided ideal generated by these line point vertices. This characterization allows us to compute the socle of some
algebras that arise as the Leavitt path algebra of some row-finite graphs. A complete description of the socle of a
Leavitt path algebra is given: it is a locally matricial algebra.
\end{abstract}

\maketitle


\section*{Introduction}

Leavitt path algebras of row-finite graphs have been recently introduced in \cite{AA1} and \cite{AMP}. They
have become a subject of significant interest, both for algebraists and for analysts working in C*-algebras.
These Leavitt path algebras $L_K(E)$ are natural generalizations of the algebras investigated by Leavitt in
\cite{L} and are a specific type of path $K$-algebras associated to a graph $E$, modulo some relations (here,
$K$ is a field).

Among the family of algebras which can be realized as the Leavitt path algebra of a graph one may find matrix rings
${\mathbb M}_n(K)$, for $n\in \mathbb{N}\cup \{\infty\}$ (where ${\mathbb M}_\infty(K)$ denotes the ring of matrices of
countable size with only a finite number of nonzero entries), the Toeplitz algebra, the Laurent polynomial ring
$K[x,x^{-1}]$, and the classical Leavitt algebras $L(1,n)$ for $n\ge 2$. Constructions like direct sums, direct limits
and matrices over the previous examples can be also achieved. We point out the reader to the papers \cite{AA1} through
\cite{APS} to get a general flavour of how to realize those algebras as the Leavitt path algebras of row-finite graphs.

In addition to the fact that these structures indeed contain many well-known algebras, one of the main interests in
their study is the comfortable pictorial representations that their corresponding graphs provide. In fact, great
efforts have been done very recently in trying to figure out the algebraic structure of $L_K(E)$ in terms of the graph
nature of $E$. Concretely, necessary and sufficient conditions on a graph $E$ have been given so that the corresponding
Leavitt path algebra $L_K(E)$ is simple \cite{AA1}, purely infinite simple \cite{AA2}, exchange \cite{APS}, finite
dimensional \cite{AAS1}, and locally finite (equivalently noetherian) \cite{AAS2}. Another approach has been the study
in \cite{AMP} of their monoids of finitely generated projective modules $V(L_K(E))$.

The socle of an algebra is a widely present notion in the mathematical literature (see \cite{D}, \cite[\S
1.1]{H}, \cite[\S IV.3]{J}, \cite[\S 7.1]{Ro}). For an algebra $A$ the (left) socle, $Soc(A)$, is defined as
the sum of all its minimal left ideals. If there are no minimal left ideals, then $Soc(A)$ is said to be
zero. When the algebra is semiprime, $Soc(A)$ coincides with the sum of all the minimal right ideals of $A$
(or it is zero in case such right ideals do not exist). It is well-known that for semiprime algebras the
socle is a sum of simple ideals; if the algebra satisfies an appropriate finiteness condition, for example
when it is left (right) artinian, then $A=Soc(A)$ is a finite direct sum of ideals each of which is a simple
left (right) artinian algebra. In this point the Wedderburn-Artin Theorem applies to describe the complete
structure of the algebra. Similar descriptions of the socle of a semiprime algebra satisfying certain chain
conditions are familiar too. Thus, if we consider the simple algebras as the building blocks, the semiprime
coinciding with their socles are the following ones.

Needless to say, despite the several steps already taken towards the understanding of the Leavitt path
algebras, no final word regarding some type of theorem of structure has been said whatsoever. In this
situation, this paper can be thought as a natural followup of the struggle of uncovering the nature of
$L_K(E)$, in the sense that a complete description of the socle of a Leavitt path algebra could lead to a
deeper knowledge of this class of algebras.

As we have already said, the Leavitt path algebras have a C*-algebra counterpart: the Cuntz-Krieger algebras
$C^*(E)$ described in \cite{R}. Both theories share many ideas and results, although they are not exactly the
same, as was revealed recently in the ``Workshop on graph algebras'' held in the University of M\'alaga (see
\cite{Work}). Because of this close connection, any advance in one field is likely to yield a breakthrough in
the other and vice versa. Thus, the results presented in this paper can be regarded as a potential tool and
source of inspiration for C*-analysts as well.

We have divided the paper into four sections. In the first one, apart from recalling some notions which will be needed
in the sequel, we show that for every graph $E$ the Leavitt path algebra $L_K(E)$ is semiprime. In sections 2 and 3 we
study the minimal left ideals of $L_K(E)$, first the ones generated by vertices (Section 2), then the general case
(Section 3).  A vertex $v$ generates a minimal left ideal if and only if there are neither bifurcations nor cycles at
any point of the tree of $v$. Such vertex $v$ will be called a line point. In general, a principal left ideal is
minimal if and only if it is isomorphic (as a left $L_K(E)$-module) to a left ideal generated by a line point.
Moreover, the set of all line points of $E$, denoted by $P_l(E)$, generates the socle of the Leavitt path algebra in
the sense that the hereditary and saturated closure of $P_l(E)$ generates $Soc(L_K(E))$ as a two-sided ideal. This is
shown in Section 4.  A complete description of the socle of a Leavitt path algebra is given: it is a locally matricial
algebra which can be seen as a Leavitt path algebra of a graph without cycles.


\section{Definitions and preliminary results}

We will first recall the graph definitions that we will need throughout the paper. For further notions on graphs we
refer the reader to \cite{AA1} and the references therein.

A \emph{(directed) graph} $E=(E^0,E^1,r,s)$ consists of two countable sets $E^0,E^1$ and maps $r,s:E^1 \to
E^0$. The elements of $E^0$ are called \emph{vertices} and the elements of $E^1$ \emph{edges}. If $s^{-1}(v)$
is a finite set for every $v\in E^0$, then the graph is called \emph{row-finite}. Throughout this paper we
will be concerned only with row-finite graphs. If $E^0$ is finite then, by the row-finite hypothesis, $E^1$
must necessarily be finite as well; in this case we say simply that $E$ is \emph{finite}. A vertex which
emits no edges (that is, which is not the source of any edge) is called a \emph{sink}. A \emph{path} $\mu$ in
a graph $E$ is a sequence of edges $\mu=e_1\dots e_n$ such that $r(e_i)=s(e_{i+1})$ for $i=1,\dots,n-1$. In
this case, $s(\mu):=s(e_1)$ is the \emph{source} of $\mu$, $r(\mu):=r(e_n)$ is the \emph{range} of $\mu$, and
$n$ is the \emph{length} of $\mu$, i.e, $l(\mu)=n$. We denote by $\mu^0$ the set of its vertices, that is:
$\mu^0=\{s(e_1),r(e_i):i=1,\dots,n\}$.

An edge $e$ is an {\it exit} for a path $\mu = e_1 \dots e_n$ if there exists $i$  such that $s(e)=s(e_i)$
and $e\neq e_i$.  If $\mu$ is a path in $E$, and if   $v=s(\mu)=r(\mu)$, then $\mu$ is called a \emph{closed
path based at $v$}. We denote by $CP_E(v)$ the set of closed paths in $E$ based at $v$.  If $s(\mu)=r(\mu)$
and $s(e_i)\neq s(e_j)$ for every $i\neq j$, then $\mu$ is called a \emph{cycle}.

For $n\ge 2$ we write $E^n$ to denote the set of paths of length $n$, and $E^*=\bigcup_{n\ge 0} E^n$ the set
of all paths. We define a relation $\ge$ on $E^0$ by setting $v\ge w$ if there is a path $\mu\in E^*$ with
$s(\mu)=v$ and $r(\mu)=w$. A subset $H$ of $E^0$ is called \emph{hereditary} if $v\ge w$ and $v\in H$ imply
$w\in H$. A hereditary set is \emph{saturated} if every vertex which feeds into $H$ and only into $H$ is
again in $H$, that is, if $s^{-1}(v)\neq \emptyset$ and $r(s^{-1}(v))\subseteq H$ imply $v\in H$. Denote by
$\mathcal{H}$ (or by $\mathcal{H}_E$ when it is necessary to emphasize the dependence on $E$) the set of
hereditary saturated subsets of $E^0$.

The set $T(v)=\{w\in E^0\mid v\ge w\}$ is the \emph{tree} of $v$, and it is the smallest hereditary subset of $E^0$
containing $v$. We extend this definition for an arbitrary set $X\subseteq E^0$ by $T(X)=\bigcup_{x\in X} T(x)$. The
\emph{hereditary saturated closure} of a set $X$ is defined as the smallest hereditary and saturated subset of $E^0$
containing $X$. It is shown in \cite{AMP} that the hereditary saturated closure of a set $X$ is
$\overline{X}=\bigcup_{n=0}^\infty \Lambda_n(X)$, where
\begin{enumerate}
\item[] $\Lambda_0(X)=T(X)$, and
\item[] $\Lambda_n(X)=\{y\in E^0\mid
s^{-1}(y)\neq \emptyset$ and $r(s^{-1}(y))\subseteq \Lambda_{n-1}(X)\}\cup \Lambda_{n-1}(X)$, for $n\ge 1$.
\end{enumerate}

We denote by $E^\infty$ the set of infinite paths $\gamma=(\gamma_n)_{n=1}^\infty$ of the graph $E$ and by $E^{\le
\infty}$ the set $E^\infty$ together with the set of finite paths in $E$ whose end vertex is a sink. We say that a
vertex $v$ in a graph $E$ is \emph{cofinal} if for every $\gamma\in E^{\le \infty}$ there is a vertex $w$ in the path
$\gamma$ such that $v\ge w$. We say that a  graph $E$ is \emph{cofinal} if so are all the vertices of $E$.

\medskip

Let $K$ be a field and $E$ a row-finite graph. We define the {\em Leavitt path $K$-algebra} $L_K(E)$ as the $K$-algebra
generated by a set $\{v\mid v\in E^0\}$ of pairwise orthogonal idempotents, together with a set of variables
$\{e,e^*\mid e\in E^1\}$, which satisfy the following relations:

(1) $s(e)e=er(e)=e$ for all $e\in E^1$.

(2) $r(e)e^*=e^*s(e)=e^*$ for all $e\in E^1$.

(3) $e^*e'=\delta _{e,e'}r(e)$ for all $e,e'\in E^1$.

(4) $v=\sum _{\{ e\in E^1\mid s(e)=v \}}ee^*$ for every $v\in E^0$ that emits edges.

In the final section of this paper many examples of Leavitt path algebras with their realizing graphs are given.
Specifically, finite (and infinite) matrix rings, matrices over classical Leavitt algebras and matrices over Laurent
polynomial algebras are built out of graphs $E$ via this $L(E)$ construction.

The elements of $E^1$ are called \emph{real edges}, while for $e\in E^1$ we call $e^\ast$ a \emph{ghost edge}.  The set
$\{e^*\mid e\in E^1\}$ will be denoted by $(E^1)^*$.  We let $r(e^*)$ denote $s(e)$, and we let $s(e^*)$ denote $r(e)$.
Unless we want to point up the base field, we will write $L(E)$ for $L_K(E)$. If $\mu = e_1 \dots e_n$ is a path, then
we denote by $\mu^*$ the element $e_n^* \dots e_1^*$ of $L (E)$.

Note that if $E$ is a finite graph then we have $\sum _{v\in E^0} v=1$; otherwise, by \cite[Lemma 1.6]{AA1}, $L (E)$ is
a ring with a set of local units consisting of sums of distinct vertices. Conversely, if $L(E)$ is unital, then $E^0$
is finite.  For any subset $H$ of $E^0$, we will denote by $I(H)$ the ideal of $L(E)$ generated by $H$.

It is shown in \cite{AA1} that $L(E)$ is a ${\mathbb Z}$-graded $K$-algebra, spanned as a $K$-vector space by $\{pq^*
\mid p,q$ are paths in $E\}$. In particular, for each $n\in \mathbb{Z}$, the degree $n$ component $L(E)_n$ is spanned
by elements of the form  $pq^*$ where $l(p)-l(q)=n$. The degree of an element $x$, denoted $deg(x)$, is the lowest
number $n$ for which $x\in \bigoplus_{m\leq n} L(E)_m$. The set of \emph{homogeneous elements} is $\bigcup_{n\in
{\mathbb Z}} L(E)_n$, and an element of $L(E)_n$ is said to be $n$-\emph{homogeneous} or \emph{homogeneous of degree}
$n$.

If $a \in L(E)$ and $d \in {\mathbb Z}^+$, then we say that $a$ is \textit{representable as an element of degree} $d$
\textit{in real (respectively ghost) edges} in case $a$ can be written as a sum of monomials from the spanning set
$\{pq^*\ \vert \ p, q \ \hbox{are paths in }E\}$, in such a way that $d$ is the maximum length of a path $p$
(respectively $q$) which appears in such monomials. Note that an element of $L(E)$ may be representable as an element
of different degrees in real (respectively ghost) edges.

The $K$-linear extension of the assignment $pq^* \mapsto qp^*$ (for $p,q$ paths in $E$) yields an involution
on $L(E)$, which we denote simply as ${}^*$. Clearly $(L(E)_n)^* = L(E)_{-n}$ for all $n\in {\mathbb Z}$.

Recall that an algebra $A$ is said to be \emph{nondegenerate} if $aL(E)a=0$ for $a\in L(E)$ implies $a=0$.

\begin{prop}\label{nondegenerate} For any graph $E$, the Leavitt path algebra $L(E)$ is nondegenerate.
\end{prop}
\begin{proof}
It is well-known that a graded algebra is nondegenerate (resp. graded nondegenerate) if and only if it is semiprime
(resp. graded semiprime). On the other hand, by \cite[Proposition II.1.4 (1)]{NvO}, a $\Z$-graded algebra is semiprime
if and only if it is graded semiprime. Hence it suffices to prove that if $a$ is any homogeneous element and
$aL(E)a=0$, then $a=0$.

For convenience we shall denote by $Z:=Z(L(E))$ the subset of elements $z\in L(E)$ such that $zL(E)z=0$. This subset
verifies $L(E)Z, ZL(E), KZ, Z^*\subseteq Z$ and does not contain neither vertices nor paths.

First we show that if $x$ is an element of $L(E)_0$, then $xL(E)x=0$ implies $x=0$. Take $0\neq x\in L(E)_0$ such that
$xL(E)x=0$ and show that this leads to a contradiction. First we analyze the trivial case in which $x$ is a linear
combination of vertices. If $v$ is one of them then $0\ne vxv\in Z$ so that we have a vertex in $Z$. Therefore $x$ is a
linear combination of vertices and of monomials $ab^*$ where $a$ and $b$ are paths of the same positive degree.

By using (4), we can always replace any vertex $w$ which is not a sink and appears in $x$, by the expression
$\sum _{\{e_i\in E^1\mid s(e_i)=w \}}e_ie_i^*.$ In that way, after simplifying if necessary, we can write $x$
as the sum of monomials of degree zero such that the only ones which are vertices are precisely sinks. In
other words, $x=x_1+x_2$, where $x_1$ is a linear combination of degree zero monomials neither of which is a
vertex, and $x_2$ is a linear combination of sinks.

Now, if we consider one of these monomials $ab^*$ appearing in the mentioned linear combination $x_1$ with
maximum degree of $a$, we can write $a=fa'$, $b=gb'$, where $f,g\in E^1$ and $a',b'$ are paths of the same
degree (in fact this degree is the degree of $a$ minus $1$).

Hence we can write $x_1=fx'g^*+z$, where $x'\in L(E)\setminus\{0\}$ and $f^*zg=0$ (this is possible because
$x_1$ contains only degree zero elements that are not vertices). Thus, by recalling that $x_2$ contains only
sinks we obtain that
$$f^*xg=f^*x_1g+f^*x_2g=f^*fx'g^*g+f^*zg+f^*x_2g=x'+0+0=x'$$ is a nonzero element of $Z$. Applying recursively to
$x'$ the argument above we get that $Z$ contains a nonzero linear combination of vertices.

To finish the proof suppose that $Z$ does not contain nonzero homogeneous elements of positive degree $<k$ and let us
prove that it does not contain nonzero homogeneous elements of degree $k$. Thus consider $0\ne x\in L(E)_k\cap Z$. For
any $f\in E^1$ we have $f^*x\in Z$ and this is an homogeneous element of degree $<k$. Therefore $f^*x=0$ for any $f\in
E^1$. Applying (4), this implies that $vx=0$ for any vertex $v$ such that $s^{-1}(v)\ne\emptyset$. On the other hand if
$v\in E^0$ is such that $s^{-1}(v)=\emptyset$, then for any $g\in E^1$ we have $vg=vs(g)g=0$ since $v\ne s(g)$. Thus
$vx=0$ for any vertex $v$ and this implies $x=0$ since $L(E)$ has local units.

Since $L(E)_{-n}=(L(E)_n)^*$, it follows that $Z$ does not contain nonzero homogeneous elements of negative degree.
\end{proof}


\section{Minimal left ideals generated by vertices}

Our first concern will be to investigate which are the conditions on a vertex $v\in E^0$ that makes minimal the left
ideal $L(E)v$. First we need the concepts of bifurcation and line point.

\begin{defis}{\rm We say that a vertex $v$ in $E^0$ is a \emph{bifurcation} (or that \emph{there is a bifurcation at} $v$)
if $s^{-1}(v)$ has at least two elements. A vertex $u$ in $E^0$ will be called a \emph{line point} if there are neither
bifurcations nor cycles at any vertex $w\in T(u)$. We will denote by $P_l(E)$ the set of all line points in $E^0$. We
say that a path $\mu$ \emph{contains no bifurcations} if the set $\mu^0\setminus\{r(\mu)\}$ contains no bifurcations,
that is, if none of the vertices of the path $\mu$, except perhaps $r(\mu)$, is a bifurcation.}
\end{defis}

\begin{lem}\label{nobifurcationsutov} Let $u,v$ be in $E^0$, with $v\in T(u)$. If the (only) path that joins $u$ with $v$
contains no bifurcations, then $L(E)u\cong L(E)v$ as left $L(E)$-modules.
\end{lem}
\begin{proof} Let $\mu\in E^*$ be such that $s(\mu)=u$ and $r(\mu)=v$. Define the right
multiplication maps $\rho_{\mu}:L(E)u\to L(E)v$ and $\rho_{\mu^*}:L(E)v\to L(E)u$, respectively, by
$\rho_{\mu}(\alpha u)=\alpha u \mu\in L(E)v$ and $\rho_{\mu^*}(\beta v)=\beta v \mu^*\in L(E)u$, for $\alpha,
\beta \in L(E) $. The fact that there are no bifurcations along the path $\mu$ allows us to apply relation
(4) to yield $\mu\mu^*=u$. Since the relation $\mu^*\mu=v$ always holds by (3), we have that $\rho_{\mu^*}
\rho_\mu =\hbox{Id}|_{L(E)u}$ and $\rho_{\mu}\rho_{\mu^*}=\hbox{Id}|_{L(E)v}$. Thus, these maps are the
desired $L(E)$-module isomorphisms.
\end{proof}

\begin{prop}\label{prior} Let $u$ be a vertex which is not a sink, and consider the set  $s^{-1}(u)=\{f_1,\ldots,f_n\}$.
Then $L(E)u=\bigoplus_{i=1}^n L(E)f_if_i^*$. Furthermore, if $r(f_i)\ne r(f_j)$ for $i\ne j$ and $v_i:=r(f_i)$, we have
$L(E)u\cong\bigoplus_{i=1}^n L(E)v_i$.
\end{prop}
\begin{proof} For $i=1,\ldots,n$, the elements $f_if_i^*$ are orthogonal idempotents by (3).  Since their sum is $u$ by
relation (4), we have $L(E)u=\bigoplus_{i=1}^n L(E)f_if_i^*$. For the second assertion in the proposition
take into account that the map $\Lambda:L(E)u\to\bigoplus_{i=1}^n L(E)v_i$ such that $x\mapsto \sum_i xf_i$
is clearly a left $L(E)$-modules homomorphism. But $\ker(\Lambda)=0$ since $\sum_i xf_i=0$ implies, by
multiplying on the right hand side by $r(f_i)$, that $xf_i=0$ for each $i$ and then $xf_if_i^*=0$. Hence
summing in $i$ we have, by relation (4), that $x=xu=\sum_i xf_if_i^*=0$. The map $\Lambda$ is also an
epimorphism since for any collection of elements $y_i\in L(E)v_i$ we have $\sum_iy_i=\Lambda(\sum_i
y_if_i^*)$.
\end{proof}

Recall that a left ideal $I$ of an algebra $A$ is said to be \emph{minimal} if it is nonzero and the only left ideals
of $A$ that it contains are $0$ and $I$. From the results above we get an immediate consequence.

\begin{corol} \label{lemmaA} Let $w$ be in $E^0$. If $T(w)$ contains some bifurcation, then the left ideal $L(E)w$
is not minimal.
\end{corol}
\begin{proof}
Let $v\in T(w)$ be a bifurcation. Consider a path $\mu=e_1\dots e_n$ joining $w$ to $v$. Take $x\in \mu^0$ the first
bifurcation occurring in $\mu$. If $x=w$ we simply apply Proposition \ref{prior}. Suppose then that $x\neq w$, so that
$x=r(e_i)$ for some $1\leq i\leq n$ and the path $\nu=e_1\dots e_i$ contains no bifurcations. Now by Lemma
\ref{nobifurcationsutov} we get $L(E)w\cong L(E)x$ as left $L(E)$-modules and by Proposition \ref{prior} we get that
$L(E)x$ is not minimal.
\end{proof}

Next we investigate another necessary condition on a vertex to generate a minimal left ideal. This is given by the
following result.

\begin{prop}\label{closedpathminimal}
If there is some closed path based at $u\in E^0$, then $L(E)u$ is not a minimal left ideal.
\end{prop}
\begin{proof} Consider $\mu\in\hbox{CP}(u)$ and suppose that $L(E)u$ is minimal. By Corollary \ref{lemmaA} there
are no bifurcations at any vertex of the path $\mu$. In particular $\mu$ is a cycle.

Consider the left ideal $0\ne L(E)(\mu+u)\subseteq L(E)u$. Since $L(E)u$ is minimal we have $u\in L(E)(\mu+u)$, so
$u=\sum_i k_i\tau_i(\mu +u)$ being each $\tau_i$ a nonzero monomial in $L(E)$ and $k_i\in K$. Note that $\tau_i\ne 0$
and $r(\tau_i)=u=s(\tau_i)$. Thus, since the tree $T(u)$ contains no bifurcations by Corollary \ref{lemmaA}, with
similar computations to that performed in \cite[Proof of Theorem 3.11]{AA1}, we see that each monomial $\tau_i$ is
either a power of $\mu$, a power of $\mu^*$ or simply $u$. Hence we have $u=p(\mu,\mu^*)(\mu+u)$, where $p$ is a
polynomial of the form
$$p(\mu,\mu^*)=l_m\mu^m+\cdots+l_1\mu+l_0u+l_{-1}\mu^*+\cdots +l_{-n}(\mu^*)^n,$$ being each $l_i$ a scalar and $m, n\ge
0$.

Taking into account that $\mu^*\mu=u=\mu\mu^*$ by relations (3) and (4), multiplying on the right by $\mu^n$ we get
$$\mu^n=(l_m\mu^{m+n}+\cdots+l_{-n}u)(\mu+u).$$ But the subalgebra of $L(E)$ generated by $\mu$ (and $u$) is isomorphic to
the polynomial algebra $K[x]$, so the previous equation implies that in $K[x]$ we have $x^n=q(x)(x+1)$ for some
polynomial $q(x)\in K[x]$. However this is impossible since evaluating in $x=-1$ we get a contradiction.
\end{proof}

Thus we have the following proposition, which gives the necessary condition on a vertex $u$ so that $L(E)u$ is a
minimal left ideal.

\begin{prop}\label{sufficient}
Let $u$ be a vertex of the graph $E$ and suppose that the left ideal $L(E)u$ is minimal. Then $u\in P_l(E)$.
\end{prop}
\begin{proof}
Take $v\in T(u)$. If there is a bifurcation at $v$ then, by Corollary \ref{lemmaA}, we get a contradiction. If there is
a cycle based at $v$, then Proposition \ref{closedpathminimal} shows that $L(E)v$ is not a minimal left ideal.
Corollary \ref{lemmaA} gives that there are no bifurcations in the (unique) path joining $u$ to $v$ so that Lemma
\ref{nobifurcationsutov} yields $L(E)u\cong L(E)v$, the former being minimal but not the latter, a contradiction.
\end{proof}

As we shall prove in what follows, this necessary condition turns out to be also sufficient.

\begin{prop}\label{corner}
For any $u\in E^0$, the left ideal $L(E)u$ is minimal if and only if $uL(E)u=K u\ \cong K$.
\end{prop}
\begin{proof}
Take into account that an element in $uL(E)u$ is a linear combination of elements of the form $k \mu$, with
$k\in K$ and $\mu$ being the trivial path $u$ or $f_1\cdots f_rg_1^*\cdots g_s^*=f_1\cdots f_r(g_s\cdots
g_1)^*$, where $f_i$ and $g_j$ are real edges and $s(f_1)=s(g_s)=u$. Apply that $T(u)$ has no bifurcations,
by Corollary \ref{lemmaA}, to obtain $f_1=g_s$, $f_2=g_{s-1}$ and so on. If $r<s$, then $\mu=f_1\dots
f_rg_s^*\dots g_{r+1}^*f_r^*\dots f_1^*$ and for $w:=r(f_r)$ we have $g_{r+1}\dots g_s\in CP(w)$. But this is
a contradiction because $w\in T(u)$ and $u\in P_l(E)$ by Proposition \ref{sufficient}. The case $r>s$ does
not happen, as can be shown analogously. Hence, $\mu=f_1\dots f_rf_r^*\dots f_1^*=u$ (there are no
bifurcations in $f_1\dots f_r$) and we have proved that $uL(E)u=Ku$.

Conversely, if $uL(E)u\cong K$, then $L(E)u$ is a minimal left ideal because for a nonzero element $au\in L(E)u$ we
have $L(E)au=L(E)u$. To show this, it suffices to prove that $u\in L(E)au$. By nondegeneracy of $L(E)$ (see Proposition
\ref{nondegenerate}), $auL(E)au\neq 0$. Take $0\neq uxau$ and apply that $uL(E)u$ is a field to obtain $ubu\in uL(E)u$
such that $u=ubuxau \in L(E)au$.
\end{proof}

\begin{rema}\label{sink}
{\rm For any sink $u$, trivially $uL(E)u=K u\cong K$, and therefore the left ideal $L(E)u$ is minimal. Also, if $w$ is
a vertex connected to a sink $u$ by a path without bifurcations, then we have that $L(E)w$ is a  minimal left ideal
because $L(E)w\cong L(E)u$ by Lemma \ref{nobifurcationsutov}.}
\end{rema}

\begin{theor}\label{lolicandido} Let $u\in E^0$. Then $L(E)u$ is a minimal left ideal if and only if $u\in P_l(E)$.
\end{theor}
\begin{proof}
Suppose that $u\in P_l(E)$. Observe that if the tree $T(u)$ is finite, then $L(E)u$ is, trivially, a minimal
left ideal, by Remark \ref{sink}, because in this case $u$ connects to a sink.

In order to prove the result for any graph $E$ we use the notion of complete subgraph given in \cite[p.
3]{AMP}. It is proved there that the row-finite graph $E$ is the union of a directed family of finite
complete subgraphs $\{E_i\}_{i\in I}$ and that the Leavitt path algebra $L(E)$ is the limit of the directed
family of Leavitt path algebras $\{L(E_i)\}_{i\in I}$ with transition monomorphisms $\varphi_{ji}:L(E_i)\to
L(E_j)$, for $i\leq j$ induced by inclusions $E_i\hookrightarrow E_j$. Denote by $\varphi_i:L(E_i)\to L(E)$
the canonical monomorphism such that $\varphi_j\varphi_{ji}=\varphi_i$ whenever $i\le j$.

To prove the minimality of $L(E)u$ we show that $uL(E)u=Ku$ and apply Proposition \ref{corner}. There is an $i\in I$
and $u_i\in L(E_i)$ such that $u=\varphi_i(u_i)$. Thus for any $a\in L(E)$ we also have $a=\varphi_j(a_j)$ for some
$j\in I$. Now, there is some $k\ge i,j$ and the tree $T(\varphi_{ki}(u_i))$ contains neither bifurcations nor closed
paths in $E_k$ since this is a subgraph of $E$. Therefore the left ideal $L(E_k)\varphi_{ki}(u_i)$ is minimal because
the graph $E_k$ is finite. Consequently $\varphi_{ki}(u_i)L(E_k)\varphi_{ki}(u_i)=K\varphi_{ki}(u_i)$ by Proposition
\ref{corner}, so that $\varphi_{ki}(u_i)\varphi_{kj}(a_j)\varphi_{ki}(u_i)=\lambda \varphi_{ki}(u_i) $ for some scalar
$\lambda\in K$. Applying $\varphi_k$ we get $uau=\lambda u$ as desired.

The converse is Proposition \ref{sufficient}.
\end{proof}

It was shown in Corollary \ref{lemmaA} that if for a vertex $u$ the tree $T(u)$ contains bifurcations, then
$L(E)u$ is not a minimal left ideal. The following example shows that the condition of not having cycles at
any point in $T(v)$ cannot be dropped in the theorem before.

\begin{exem} \label{aasdosexample} {\rm Consider the graph $E$ given by $$\xymatrix{{\bullet}^u \ar[r]^e & {\bullet}^v
\ar@(dr,ur)_f}$$ Then $L(E)u$ is {\it not} a minimal left ideal (note that  there {\it is} a cycle in $v\in
T(u)$). To show this we use \cite[Theorem 3.3]{AAS2} to get that $L(E)\cong A:={\mathbb M}_2(K[x,x^{-1}])$
via an isomorphism which sends $u$ to $e_{22}=\left(\begin{smallmatrix} 0 & 0 \\ 0 & 1
\end{smallmatrix}\right)\in A$. Now if $L(E)u$ were a minimal left ideal, then so would be $Ae_{22}$, but the
nonzero left ideal (of $A$) $I=\left(\begin{smallmatrix}0& \langle 1+x\rangle\\ 0 & \langle
1+x\rangle\end{smallmatrix}\right)$
is strictly contained in $Ae_{22}= \left(\begin{smallmatrix} 0 & K[x,x^{-1}] \\
0 & K[x,x^{-1}] \end{smallmatrix}\right)$, a contradiction.}
\end{exem}


\section{Minimal left ideals}

The following result is the key tool to obtain the reduction process needed to translate the minimality of a principal
left ideal to a left ideal generated by a vertex. Moreover, it can be used to shorten the proof given in \cite{AA1} to
show that if a graph $E$ satisfies \emph{Condition} (L) (that is, if every cycle has an exit) and the only hereditary
and saturated subsets of $E^0$ are the trivial ones, then the associated Leavitt path algebra is simple.

\begin{prop}\label{cases} Let $E$ be a graph. For every nonzero element $x\in L(E)$ there exist
$\mu_1,\dots,\mu_r,\nu_1,\dots,\nu_s\in E^0\cup E^1\cup (E^1)^*$ such that:
\begin{enumerate}
\item $\mu_1\dots\mu_rx\nu_1\dots\nu_s$ is a nonzero element in $Kv$, for some $v\in E^0$, or
\item  there exist a vertex $w$ and a cycle without exits $c$ based at $w$ such that
$\mu_1\dots\mu_rx\nu_1\dots\nu_s$ is a nonzero element in $wL(E)w=\{\sum_{i=-m}^n k_i c^i\mbox{ for }m,n\in
{\mathbb N } \mbox{ and } k_i \in K\}$.
\end{enumerate}
Both cases are not mutually exclusive.
\end{prop}
\begin{proof}
Show first that for a nonzero element $x \in L(E)$, there exists a path $\mu\in L(E)$ such that $x\mu$ is nonzero and
in only real edges.

Consider a vertex $v\in E^0$ such that $x v\neq 0$. Write $x v= \sum_{i=1}^m \beta_ie_i^\ast +\beta$, with $e_i\in
E^1$, $e_i\neq e_j$ for $i\neq j$ and $\beta_i, \beta \in L(E)$, $\beta$ in only real edges and such that this is a
minimal representation of $x v$ in ghost edges.

If $x ve_i=0$ for every $i\in \{1, \dots, m\}$, then $0=x ve_i= \beta_i +\beta e_i$, hence $\beta_i=-\beta e_i$, and $x
v=\sum_{i=1}^m -\beta e_i e_i^\ast +\beta=\beta (\sum_{i=1}^m - e_i e_i^\ast +v)\neq 0$. This implies that
$\sum_{i=1}^m - e_i e_i^\ast +v\neq 0$ and since $s(e_i)=v$ for every $i$, this means that there exists $f\in E^1$,
$f\neq e_i$ for every $i$, with $s(f)=v$. In this case, $x v f =\beta f \neq 0$ (because $\beta$ is in only real
edges), with $\beta f$ in only real edges, which would conclude our discussion.

If $x ve_i\neq 0$ for some $i$, say for $i=1$, then $0\neq x v e_1=\beta_1+\beta e_1$, with $\beta_1+\beta e_1$ having
strictly less degree in ghost edges than $x$.

Repeating this argument, in a finite number of steps we prove our first statement.

\medskip

Now, assume  $x=xv$ for some $v\in E^0$ and $x$ in only real edges. Let $0\ne x=\sum_{i=1}^r k_i\alpha_i$ be
a linear combination of different paths $\alpha_i$ with $k_i\ne 0$ for any $i$. We prove by induction on $r$
that after multiplication on the left and/or the right we get a vertex or a polynomial in a cycle with no
exit. For $r=1$ if $\alpha_1$ has degree $0$ then it is a vertex and we have finished. Otherwise we have
$x=k_1\alpha_1=k_1 f_1\cdots f_n$ so that $k_1^{-1}f_n^*\cdots f_1^*x=v$ where $v=r(f_n)\in E^0$.

Suppose now that the property is true for any nonzero element which is a sum of less than $r$ paths in the
conditions above. Let $0\ne x=\sum_{i=1}^r k_i\alpha_i$ such that $\deg(\alpha_i)\le\deg(\alpha_{i+1})$ for
any $i$.  If for some $i$ we have $\deg( \alpha_i)=\deg(\alpha_{i+1})$ then, since $\alpha_i\ne\alpha_{i+1}$,
there is some path $\mu$ such that $\alpha_i=\mu f\nu$ and $\alpha_{i+1}=\mu f'\nu'$ where $f,f'\in E^1$ are
different and $\nu,\nu'$ are paths. Thus $0\ne f^*\mu^*x$ and we can apply the induction hypothesis to this
element. So we can go on supposing that $\deg(\alpha_i)<\deg(\alpha_{i+1})$ for each $i$.

We have  $0\ne \alpha_1^*x=k_1v+\sum_i k_i\beta_i$, where $v=r(\alpha_1)$ and $\beta_i=\alpha_1^*\alpha_i$. If some $\beta_i$ is null then apply the
induction hypothesis to $\alpha_1^*x$ and we are done. Otherwise if some $\beta_i$ does not start (or finish) in $v$ we apply the induction
hypothesis to $v\alpha_1^*x\ne 0$ (or $\alpha_1^*xv\ne 0$). Thus we have $$0\ne z:=\alpha_1^* x=k_1v+\sum_{i=1}^r k_i\beta_i,$$ where
$0<\deg(\beta_1)<\cdots< \deg(\beta_r)$ and all the paths $\beta_i$ start and finish in $v$.

Now, if there is a path $\tau$ such that $\tau^*\beta_i=0$ for some $\beta_i$ but not for all of them, then we apply
our inductive hypothesis to $0\ne \tau^*z\tau$. Otherwise for any path $\tau$ such that $\tau^*\beta_j=0$ for some
$\beta_j$, we have $\tau^*\beta_i=0$ for all $\beta_i$. Thus $\beta_{i+1}=\beta_i r_i$ for some path $r_i$ and $z$ can
be written as
$$z=k_1v+k_2\gamma_1+k_3\gamma_1\gamma_2+\cdots+k_r\gamma_1\cdots
\gamma_{r-1},$$ where each path $\gamma_i$ starts and finishes in $v$. If the paths $\gamma_i$ are not
identical we have $\gamma_1\ne\gamma_i$ for some $i$, then $0\ne\gamma_i^*z\gamma_i=k_1v$ proving our thesis.
If the paths are identical then $z$ is a polynomial in the cycle $c=\gamma_1$ with independent term $k_1v$,
that is, an element in $vL(E)v$.

If the cycle has an exit, it can be proved that there is a path $\eta$ such that $\eta^*c=0$, in the following way: Suppose that there is a vertex
$w\in T(v)$, and two edges $e, f$, with $e\neq f$, $s(e)=s(f)=w$,  and such that $c=aweb=aeb$, for $a$ and $b$ paths in $L(E)$. Then $\eta=af$ gives
$\eta^*c=f^*a^*aeb=f^*eb=0$. Therefore, $\eta^*z\eta$ is a nonzero scalar multiple of a vertex.

Moreover, if $c$ is a cycle without exits, with similar ideas to that of \cite[Proof of Theorem 3.11]{AA1}, it is not
difficult to show that
$$vL(E)v=\left\{\sum_{i=-m}^nl_ic^i, \ \hbox{with }l_i\in K \ \hbox{and }  m, n \in \mathbb{N}\right\},$$ where we
understand $c^{-m}=(c^*)^m$ for $m\in \mathbb{N}$ and $c^0=v$.

Finally, consider the graph $E$ consisting of one vertex and one loop based at the vertex to see that both cases can
happen at the same time. This completes the proof.
\end{proof}

\begin{corol}Let $E$ be a graph that satisfies Condition {\rm (L)} and such that the only hereditary and saturated
subsets of $E^0$ are the trivial ones. Then the associated Leavitt path algebra is simple.
\end{corol}
\begin{proof} Let $I$ be a nonzero ideal of $L(E)$. By Proposition \ref{cases}, $I\cap E^0\neq \emptyset$.
Since $I\cap E^0$ is hereditary and saturated (\cite[Lemma 3.9]{AA1}), it coincides with $E^0$. This means
$I=L(E)$.
\end{proof}

The following result plays an important role in the proof of the main result of \cite{AA2}, that
characterizes those graphs $E$ for which the Leavitt path algebra is purely infinite and simple (see
\cite[Proposition 6]{AA2}).

\begin{corol}\label{caseswithL} If a graph $E$ satisfies Condition {\rm (L)}, then for every nonzero element $x\in L(E)$
there exist $\mu_1,\dots,\mu_r,\nu_1,\dots,\nu_s\in E^0\cup E^1\cup (E^1)^*$ and $v\in E^0$ such that $0 \neq
\mu_1\dots\mu_rx\nu_1\dots\nu_s\in Kv$.
\end{corol}

\begin{theor}\label{reduction} Let $x$ be in $L(E)$ such that $L(E)x$ is a minimal left ideal. Then,
there exists a vertex $v\in P_l(E)$ such that $L(E)x$ is isomorphic (as a left $L(E)$-module) to $L(E)v$.
\end{theor}
\begin{proof}
Consider $x\in L(E)$ as in the statement. By Proposition \ref{cases} we have two cases. Let us prove that the second
one is not possible.

Suppose, otherwise, that there exist a vertex $w$ and a cycle without exits $c$ based at $w$ such that
$\lambda:=\mu_1\dots\mu_rx\nu_1\dots\nu_s\in wL(E)w=\{\sum_{i=-m}^n k_i c^i\mbox{ for some }m,n\in {\mathbb N},
\hbox{and } k_i \in K\}$. Note that $wL(E)w$ is isomorphic to $K[t, t^{-1}]$ as a $K$-algebra and that $\varphi: K[t,
t^{-1}]\to L(E)$ given by $\varphi(1)=w$, $\varphi(t)=c$ and $\varphi(t^{-1})=c^*$, is a monomorphism with image
$wL(E)w$. Since $L(E)\lambda $ is isomorphic to $L(E)x$, then it is a minimal left ideal of $L(E)$. (Note that
$L(E)x=L(E)\mu_1\dots\mu_rx$ by the minimality of $L(E)x$; moreover, for $\nu:=\nu_1\dots\nu_s$, the map $\rho_\nu:
L(E)x \to L(E)x\nu$ given by $\rho_\nu(y)=y\nu$ is a nonzero epimorphism of left $L(E)$-modules. The simplicity of
$L(E)x$ implies that it is an isomorphism.) Now, consider $wL(E)\lambda$, which is a minimal left ideal of $wL(E)w$.
Then the nonzero left ideal $\varphi^{-1}(wL(E)\lambda)$ is minimal in $K[t, t^{-1}]$, a contradiction, since this
algebra has no minimal left ideals.

Hence, we are under case (1) of Proposition \ref{cases}, and so there exist $\mu_1,\dots,\mu_r,\nu_1,\dots,\nu_s\in
E^0\cup E^1\cup (E^1)^*$, $k\in K$, such that $0\neq \mu_1\dots\mu_rx\nu_1\dots\nu_s = kv$, for some $v\in E^0$. Then
$L(E)v=L(E)kv=L(E)\mu_1\dots\mu_rx\nu_1\dots\nu_s\cong L(E)x$, as left $L(E)$-modules, as required. Finally, apply
Theorem \ref{lolicandido} to obtain that $v\in P_l(E)$.
\end{proof}


\section{The socle of a Leavitt path algebra}

Having characterized in the previous section the minimal left ideals, we are in a position to finally compute, in this
section, the socle of a Leavitt path algebra. We will achieve this by giving a generating set of vertices of the socle
as a two-sided ideal.

\begin{prop}\label{leftsocle} For a graph $E$ we have that $\sum_{u\in P_l(E)}L(E)u\subseteq Soc(L(E))$. The
reverse containment does not hold in general.
\end{prop}
\begin{proof} By Theorem \ref{lolicandido}, given $u\in P_l(E)$, the left ideal $L(E)u$ is minimal and therefore it is
contained in the socle.

We exhibit an example to show that the converse containment is not true: consider the graph $E$ given by
$$\xymatrix{ {\bullet}^v & {\bullet}^z \ar[l]_e \ar[r]^f & {\bullet}^w }$$ By \cite[Proposition 3.5]{AAS1},
the Leavitt path algebra of this graph is $L(E)\cong {\mathbb M}_2(K)\oplus {\mathbb M}_2(K)$, and therefore
it coincides with its socle. However, $Soc(L(E))=L(E)\neq \sum_{u\in P_l(E)}L(E)u=L(E)v+L(E)w$ as for
instance $e^*\not\in L(E)v+L(E)w$. (To see this, suppose that $e^*=\alpha v + \beta w$, then $e^*=e^*z=\alpha
vz + \beta wz=0$, a contradiction.)
\end{proof}

Nevertheless, although the previous result shows that in general the socle of a Leavitt path algebra is not necessarily
the principal left ideal generated by $P_l(E)$, it turns out that the socle of a Leavitt path algebra is indeed the
two-sided ideal generated by this set of line points $P_l(E)$.

\begin{theor}\label{socleofLPA} Let $E$ be a graph. Then $Soc(L(E))=I(P_l(E))=I(H)$, where $H$ is the hereditary and
saturated closure of $P_l(E)$.
\end{theor}
\begin{proof} First we show that $Soc(L(E))=I(P_l(E))$. Take a minimal left ideal $I$ of $L(E)$. The Leavitt path algebra
$L(E)$ is nondegenerate  (Proposition \ref{nondegenerate}), therefore an standard argument shows that there exists
$\alpha=\alpha^2 \in L(E)$ (not necessarily a vertex) such that $I=L(E)\alpha$.

Apply Proposition \ref{reduction} to get that $L(E)\alpha\cong L(E)u$ for some $u\in P_l(E)$. Let $\phi:
L(E)\alpha\to L(E)u$ be an $L(E)$-module isomorphism. Write $\phi(\alpha)=xu$ and $\phi^{-1}(u)=y\alpha$ for
some $x,y\in L(E)$; thus: $\alpha=\phi^{-1}\phi(\alpha)=\phi^{-1}(xu^2)=xu\phi^{-1}(u)=xuy\alpha$.
Analogously we have $u=y\alpha xu$. Then, by naming $a=xu$ and $b=y\alpha$, we get that $\alpha=ab$ and
$u=ba$, for some $a,b\in L(E)$. Hence, $\alpha=abab=aub\in I(P_l(E))$.

To see the converse containment pick $v\in P_l(E)$ and show that $L(E)vL(E)\subseteq Soc(L(E))$. By Proposition
\ref{leftsocle} we have that $L(E)v\subseteq Soc(L(E))$; since the socle is always a two-sided ideal, we have our
claim.

Finally, apply \cite[Lemma 2.1]{APS} to obtain that $I(P_l(E))=I(\overline{P_l(E)})$, where $H=\overline{P_l(E)}$ is
indeed the hereditary and saturated closure of $P_l(E)$.
\end{proof}

This result has an immediate but useful corollary.

\begin{corol}\label{nonzerosocle} For a graph $E$, the Leavitt path algebra $L(E)$ has nonzero socle if and
only if $P_l(E)\neq \emptyset$.
\end{corol}

We obtain some consequences of this result. The first one is that arbitrary matrix rings over the classical
Leavitt algebras $L(1,n)$, for $n\geq 2$, as well as over the Laurent polynomial algebras $K[x,x^{-1}]$, all
have zero socle. The second is that for Leavitt path algebras of finite graphs (this class in particular
includes the locally finite, or equivalently, noetherian, Leavitt path algebras studied in \cite{AAS2}) we
can find a more specific necessary and sufficient condition so that they have nonzero socle.

\begin{corol}\label{classicalLeavitt} For all $m,n\geq 1$, $Soc({\mathbb M}_m(L(1,n)))=0$.
\end{corol}
\begin{proof} By taking into account both \cite[Proposition 12]{AA2} for the case $n\geq 2$ and \cite[Theorem
3.3]{AAS2} for the case $n=1$, the know that the algebra $A={\mathbb M}_m(L(1,n))$ is the Leavitt path algebra of the
graph $E_n^m$ given by $$
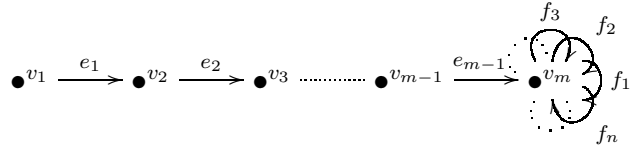
\xymatrix{{\bullet}^{v_1} \ar [r] ^{e_1} & {\bullet}^{v_2} \ar [r] ^{e_2} & {\bullet}^{v_3}
\ar@{.}[r] & {\bullet}^{v_{m-1}} \ar [r] ^{e_{m-1}} & {\bullet}^{v_m}
 \ar@(ur,dr) ^{f_1} \ar@(u,r) ^{f_2} \ar@(ul,ur) ^{f_3} \ar@{.} @(l,u) \ar@{.} @(dr,dl) \ar@(r,d) ^{f_n}& }$$ This
graph clearly has $P_l(E_n^m)=\emptyset$, so that Corollary \ref{nonzerosocle} gives the result.
\end{proof}

\begin{corol}\label{locallyfinitesocle} Let $L(E)$ be a Leavitt path algebra with $E$ a finite graph. Then $L(E)$
has nonzero socle if and only if $E^0$ has a sink.
\end{corol}
\begin{proof} If $L(E)$ has nonzero socle, Corollary \ref{nonzerosocle} gives that $P_l(E)\neq \emptyset$.
Take $v\in P_l(E)$. Since $T(v)$ has no bifurcations, contains no cycles and the graph is finite, clearly $T(v)$ must
contain a sink. Conversely, any sink $w$ obviously has $w\in P_l(E)$, so that Corollary \ref{nonzerosocle} gives
$Soc(L(E))\neq 0$.
\end{proof}

It is well-known that for $A_n:={\mathbb M}_n(K)$, with $n\in {\mathbb N}\cup \{\infty\}$, then $A_n$ coincides with
its socle. Theorem \ref{socleofLPA} can be applied to obtain these results by using the Leavitt path algebra approach.
Concretely, if $n$ is finite then $A_n$ is the Leavitt path algebra of the finite line graph $E_n$ given by
$$\xymatrix{ {\bullet}^{v_1} \ar[r] & {\bullet}^{v_2} \ar@{.}[r] & {\bullet}^{v_{n-1}}\ar[r]  & {\bullet}^{v_n}}$$
Whereas $A_\infty$ can be realized as $L(E_\infty)$ for the infinite graph $E_\infty$ defined as
$$\xymatrix{ {\bullet}^{v_1} \ar[r] & {\bullet}^{v_2} \ar[r] & {\bullet}^{v_3} \ar@{.}[r]  & }$$
In any case, clearly $P_l(E_n)=E_n^0$, so that Theorem \ref{socleofLPA} applies to give
$Soc(A_n)=I(E_n^0)=L(E_n)=A_n$, since the sum of vertices is a set of local units for $L(E_n)$.

We can perform analogous computations with arbitrary algebras of the form $\bigoplus_{i\in I} {\mathbb M}_{n_i}(K)$, where $I$ is any countable set
and $n_i\in {\mathbb N}\cup \{\infty\}$ for every $i\in I$ since these can be realized as the Leavitt path algebras of disjoint unions of graphs of
the form above, for which all its vertices are line points.

\begin{exem} {\rm Not every acyclic graph coincides with its socle. Let $E$ be the following graph:}

$$
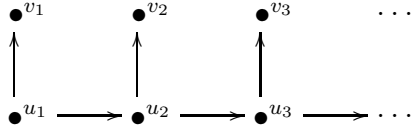
\xymatrix{
\bullet^{v_1} & \bullet^{v_2} & \bullet^{v_3}  & \dots & \\
\ar@<1ex>[u] \bullet^{u_1} \ar[r] & \ar@<1ex>[u]\bullet^{u_2}  \ar[r]& \ar[r] \ar@<1ex>[u]\bullet^{u_3}&
\ar@{.}\dots}$$

 {\rm We claim that $L(E)$ does not coincide with its socle. Otherwise, by Theorem \ref{socleofLPA},
$L(E)=I(H)$, where $H$ is the hereditary and saturated closure of $P_l(E)=\{v_n \ \vert\ n \in \mathbb{N}\}$. It is not
difficult to see that $P_l(E)$ is hereditary and saturated, hence $H=P_l(E)$. By \cite[Theorem 4.3]{AMP} $I(H)=I(E^0)$
implies $H=E^0$, a contradiction.}
\end{exem}

We finish the paper giving a complete characterization of the socle of a Leavitt path algebra.

Recall that a \emph{matricial algebra} is a finite direct product of full matrix algebras over $K$, while a
\emph{locally matricial algebra} is  a direct limit of matricial algebras.

The following definitions are particular cases  of those appearing in \cite[Definition 1.3]{DHSz}:

Let $E$ be a graph, and let $\emptyset \ne H\in \mathcal{H}_E$. Define
$$F_E(H)=\{ \alpha =(\alpha_1, \dots ,\alpha_n)\mid \alpha _i\in
E^1, s(\alpha _1)\in E^0\setminus H, r(\alpha _i)\in E^0\setminus H \mbox{ for } i<n, r(\alpha _n)\in H\}.$$
Denote by $\overline{F}_E(H)$ another copy of $F_E(H)$. For $\alpha\in F_E(H)$, we write $\overline{\alpha}$
to denote a copy of $\alpha$ in $\overline{F}_E(H)$. Then, we define the graph ${}_HE=({}_HE^0, {}_HE^1, s',
r')$ as follows:
\begin{enumerate}
\item $({}_HE)^0=H\cup F_E(H)$. \item $({}_HE)^1=\{ e\in
E^1\mid s(e)\in H\}\cup \overline{F}_E(H)$.
\item For every $e\in E^1
\mbox{ with } s(e)\in H$, $s'(e)=s(e)$ and $r'(e)=r(e)$. \item For every $\overline{\alpha}\in
\overline{F}_E(H)$, $s'(\overline{\alpha})=\alpha$ and $r'(\overline{\alpha})=r(\alpha)$.
\end{enumerate}

\begin{theor}\label{estructuradelzocalo}
For a graph $E$ the socle of the Leavitt path algebra $L(E)$ is a locally matricial algebra.
\end{theor}
\begin{proof} Suppose that our graph $E$ has line points (otherwise the socle of $L(E)$ would be 0 and the
result would follow trivially). We have proved in Theorem \ref{socleofLPA} that $Soc(L(E))=I(H)$, where $H$ is the
hereditary and saturated closure of $P_l(E)$. By \cite[Lemma 1.2]{AP}, $I(H)\cong L(_HE)$. If we had proved that $_HE$
is an acyclic graph then, by \cite[Corollary 3.6]{APS}, the Leavitt path algebra $L(_HE)$ would be locally matricial,
and the proof would be complete. Hence, let us prove this statement. Suppose on the contrary that there exists a cycle
$C$ in $_HE$. By the definition of $_HE$ we have that $C$ has to be a cycle in $E$ with vertices in $H$. Let $n$ be the
smallest non-negative integer having $\Lambda_n(P_l(E))\cap C^0\neq\emptyset$. Choose $v$ in this intersection. If
$n>0$ then $\Lambda_{n-1}(P_l(E))\cap C^0=\emptyset$ and, therefore, $\emptyset\neq r(s^{-1}(v))\subseteq
\Lambda_{n-1}(P_l(E))$. In particular $\Lambda_{n-1}(P_l(E))\cap C^0\neq\emptyset$, a contradiction, so $n$ must be
zero and consequently $T(P_l(E))\cap C^0= P_l(E)\cap C^0\neq\emptyset $. But this is a contradiction because of the
definition of $P_l(E)$.
\end{proof}


\section*{acknowledgments}
The first author gratefully acknowledges the support from the Centre de Recerca Matem\`{a}tica. All authors were partially supported by the Spanish
MEC and Fondos FEDER jointly through projects MTM2004-06580-C02-02 and by the Junta de Andaluc\'{\i}a PAI projects FQM-336 and FQM-1215. Parts of
this work were carried out during visits of the first author to the Universidad de M\'alaga, Spain and to the Centre de Recerca Matem\`{a}tica,
Spain. The first author wishes to thank both host centers for their warm hospitality.


\end{document}